\documentclass[11pt,twoside]{article}
\usepackage{amsfonts,amsmath,amsthm,amssymb,cases,mathrsfs} \marginparwidth 100pt

\def\thebibliography#1{\center{\bf\normalsize References}\list
 {[\arabic{enumi}]}{\settowidth\labelwidth{[#1]}\leftmargin\labelwidth
 \advance\leftmargin\labelsep
 \usecounter{enumi}}
 \def\newblock{\hskip .11em plus .33em minus .07em}
 \sloppy\clubpenalty4000\widowpenalty4000
 \sfcode`\.=1000\relax}

\def\nnb{\nonumber}
\def\ds{\displaystyle}
\def\cd{\cdot}
\def\cds{\cdots}
\newcommand{\refeq}[1]{$(\ref{#1})$}
\newcommand{\thb}[1]{{\rm (#1)}}
\def\endpf{\hfill$\Box$\vspace{0.4cm}}

\def\eqon{ \, {\rm on } \, \,}
\def\eqin{ \, {\rm in } \, \,}
\def\eqae{ \, {\rm a.e. } \,\, }

\def\all{  \, \forall \, }
\newcommand{\si}[1]{\mbox{ strongly in}\, #1}
\newcommand{\wi}[1]{\mbox{ weakly  in}\, #1}

\def\ol{\overline}
\def\a{\alpha}
\def\b{\beta}

\def\d{\delta}

\def\ve{\varepsilon}
\def\z{\zeta}

\def\l{\lambda}

\def\m{\mu}

\def\s{\sigma}
\def\t{\tau}

\def\o{\omega}

\def\L{\Lambda}

\def\F{\Phi}
\def\O{\Omega}
\def\cA{{\cal A}}

\def\cE{{\cal E}}

\def\cG{{\cal G}}

\def\cM{{\cal M}}

\def\cP{{\cal P}}

\def\cR{{\cal R}}
\def\cS{{\cal S}}

\def\cU{{\cal U}}

\def\mcM{{\mathscr M}}

\def\barf{\bar f}

\def\bu{\bar u}

\def\by{\bar y}
\def\bz{\bar z}

\def\bGO{\ol{\O}}
\def\hA{\widehat A}

\def\tiA{\widetilde A}

\def\tiE{\widetilde E}

\def\tiw{\tilde w}
\def\tix{\tilde x}
\def\tiy{\tilde y}

\def\tiGo{\tilde \o}
\def\qq{\qquad}
\def\q{\quad}
\def\pa{\partial}
\def\na{\nabla}
\def\co{\,{\rm co\,}}
\def\coh{\,{\rm \overline{\, co}}\,}

\def\diam{\,{\rm diam}\,}

\def\limsup{\mathop{\overline{\rm lim}}}
\def\liminf{\mathop{\underline{\rm lim}}}

\def\Hto{\xrightarrow{H}}
\def\wto{\rightharpoonup}
\def\Rto{\xrightarrow{\cR}}
\def\IR{{\rm I\hspace{-0.90mm}R}\mbox{}}

\def\defeq{\buildrel \triangle \over =}
\def\dbR{{\mathop{\rm l\negthinspace R}}}

\newcommand{\beq}[2]{\begin{equation}\label{#1}#2\end{equation}}

\newcommand{\set}[1]{\left\{#1\right\}}

\newtheorem{lemma}{Lemma}[section]
\newtheorem{definition}[lemma]{Definition}
\newtheorem{theorem}[lemma]{Theorem}
\newtheorem{prop}[lemma]{Proposition}
\newtheorem{Remark}{Remark}[section]

\numberwithin{equation}{section} 

\title{Cesari-type  Conditions for Semilinear Elliptic Equations with Leading Term Containing Controls
\thanks{This work was supported in part by NSFC (No.
61074047 and 10831007), and 973 Program (No. 2011CB808002)}}

\author{\rm Bo Li\footnote{School of Mathematical Sciences, Fudan University, Shanghai 200433,
China (Email: \texttt{062018053@fudan.edu.cn})}~~~and~~~Hongwei
Lou\footnote{School of Mathematical Sciences, and LMNS, Fudan
University, Shanghai 200433, China (Email:
\texttt{hwlou@fudan.edu.cn})}}
\date{}
\allowdisplaybreaks
\begin{document}
\maketitle

\footnotesize \textbf{Abstract.} An optimal control problem governed
by semilinear elliptic partial differential equations is considered.
The equation is in divergence form with the leading term containing
controls. By studying the $G$-closure of the leading term, an
existence result is established under a  Cesari-type  condition.



\normalsize

\section{Introduction}
Consider the following controlled elliptic partial differential equation of divergence form:
\beq{orgin}{\left\{
    \begin{array}{ll}
    \ds -\na\cd(A(x,u(x))\na y(x))=f(x,y(x),u(x)), & \eqin \,\O,\\
   y(x)=0, & \eqon \,\pa\O,
    \end{array}
    \right.
}
where $\O$ is a smooth bounded domain in $\IR^n$, $A: \O\times U\to
\IR^{n\times n}$ is a map taking values in the set of all positive
definite matrices, and $f:\O\times \IR\times U\to \IR$, with $U$
being a separable metric space. The control function $u(\cd)$ is
taken from the set
$$
\cU\equiv \set{v:\O \to U \big|v(
\cd) \,\text{\,is measurable}}.
$$
Let the cost functional be defined by
 \beq{J}{J(u(\cd))=\int_\O f^0(x,y(x),u(x))\,dx,}
where $y(\cd)$ is the solution of \refeq{orgin} (called the state
corresponding to control $u(\cd)$). Our optimal control problem is
as follows.

\textbf{Problem (C)}. Find a $\bu(\cd)\in\cU$ such that
\begin{equation}\label{C}
J(\bu(\cd))=\inf_{u(\cd)\in\cU} J(u(\cd)).
\end{equation}

Any $\bu(\cd)$ satisfying \refeq{C} is called an optimal control. It
is well-known that optimal control of Problem (C) may fail to exist.
When $A(x,u)\equiv A(x)$, a suitable Cesari-type condition and some
other mild conditions will guarantee the existence of an optimal
control. Cesari-type condition  is a natural generalization of
optimal control problem with linear state equations and convex cost
functionals. Many results are available along these lines. For
further detail, see the books by Cesari \cite{Cesari}, Li and Yong
\cite{Li}, for examples. For the two phrase case, i.e.,
$U=\set{0,1}$ and $A(x,i)\equiv A_i$ $(i=0,1)$ with $A_0,A_1$ being
two constant matrices, Murat and Tartar gave an existence result in
the framework of ``relaxation" control (see \cite{Murat2}). However,
it seems no work devoted to the existence of optimal controls for
general cases. 

In this paper, we will give a Cesari-type result to ensure the
existence of a solution to Problem  (C). We always assume $\L$ and
$\l$ be two constants satisfying $\L\geq \l>0$. Denote by $\cS^n_+$
the set of all $n\times n$ (symmetric) positive definite matrices
and
$$
\mcM_{\L,\l}=\set{Q\in \cS^n_+\Big|\l |\xi|^2\leq Q \xi\cd\xi\leq
    \L|\xi|^2,\q \all\, \xi\in \IR^n}.
$$
For a matrix $B$, we always denote $\ds B_{ij}$ as its entries.

We recall that a Polish space is a  separable completely metrizable
topological space. We mention that all (nonempty) closed sets and
open sets in $\dbR^m$ are polish spaces.

We make the following assumptions.

(S1) Set $\O$ is a bounded domain in $\IR^n$ with a $C^2$ boundary
$\pa\O$.

(S2) $U$ is a Polish space.

(S3) Function  $A(x,v)$ takes values in  $\mcM_{\L,\l}$, which are
measurable in $x\in \O$ and continuous in $v\in U$. Further, there
exists an $F\in L^\infty(\O;\IR^m)$ and a continuous $\o:
 [0,+\infty)\to [0,+\infty)$, such that $\o(0)=0$ and
\beq{S3}{\big|A(x,v)-A(\tix,v)\big|\leq
    \o\Big(\big|F(x)-F(\tix)\big|\Big),\qq\all x,\tix\in\O,\,v\in U.}

(S4) Function $f(x,y,v)$ is measurable in~$x$ and continuous in
$(y,v)\in \IR\times U$ for almost all $x\in \O$. Moreover, for
almost all $x\in \O$,
\beq{fy}{f_y(x,y,v)\leq 0,\qq\all  \, (y,v)\in \IR\times U,}
and for any $R>0$, there exists an $M_R>0$ such that
 \beq{f}{|f(x,y,v)|+|f_y(x,y,v)|\leq
M_R,\qq\all v\in  U,\,|y|\leq R.}

(S5) Function $f^0(x,y,v)$ is measurable in~$x$, lower
semicontinuous in $(y,v)\in\IR\times U$ for almost all $x\in \O$.
 Moreover, for almost all $x\in \O$ and for any $R>0$, there exists an $K_R>0$ such that
\beq{}{f^0(x,y,v)\geq -K_R,\qq\all\,v\in  U,\,|y|\leq R.}
\begin{Remark} In \refeq{S3}, $m$ need not equal to $n$. We can see
that \refeq{S3} holds naturally when $U$ is a finite set. On the
other hand, if $A(x,u)$ is uniformly continuous in $x\in \O$ with
respect to $u\in U$, then \refeq{S3} holds.
\end{Remark}
\begin{Remark}\label{R102} Without loss of generality, we can suppose that $\o(\cd)$ is a
continuous module in \thb{S3}, i.e.,

\thb{i} $\o(\cd)$  is continuous and increasing on $[0,+\infty)$,

\thb{ii} $\o(0)=0$,

\thb{iii} it holds that
$$
\o(r+s)\leq \o(r)+\o(s), \qq\all r,s>0.
$$

In fact, if necessary, we can replacing $\o(\cd)$ by
$$
\tiGo(r)\defeq \sup_{|s-t|\leq r\atop s,t\in [0,R]}
\big|\o(s)-\o(t)\big|, \qq r\geq 0,
$$
where $\ds R=2\|F\|_{L^\infty(\O;\IR^m)}$.
\end{Remark}

Denote $Z=[0,1]^n$ and
$$
\cU_Z\equiv \set{v:Z \to U \big|v(\cd) \,\text{\, is
measurable}}.
$$
Let $e_1,e_2,\ldots,e_n$ be the canonical basis of $\IR^n$. We call
a function $g(x)$ is $Z$-periodic if it admits periodic $1$ in the
direction $e_j$ ($j=1,2,\ldots,n$). Denote
\begin{align*}
L^\infty_\#(Z)=\set{h\in L^\infty(\IR^n)\big|h \text{\, is $Z$-periodic}},\\
H^1_\#(Z)=\set{h\in H^1_{loc}(\IR^n)\big| h \text{\, is $Z$-periodic}}.
\end{align*}
We define
\begin{equation}\label{E}
  \cE(x,y)=\Big\{(P,\zeta, \zeta^0)\in\cS^n_+ \times \IR \times \IR \, \Big|P=A(x,u),
       \,\z=f(x,y,u),\, \zeta^0\geq f^0(x,y,u),\,u\in U \Big\}
\end{equation}
and
\begin{equation}\label{GE}
\begin{split}
   G\cE(x,y)=\Big\{(&P,\z,\, \zeta^0)\in  \cS^n_+  \times \IR \times \IR \, \Big|~P_{ij}=\int_Z A(x,u(z))(e_i+\na w^i (z;x))\cd e_j\,dz,\\
   &\text{where}~~w^i(\cd)\in H^1_\#(Z)\text{\, solves}
    ~~ \na_z\cd\Big(A(x,u(z))(e_i+\na w^i (z;x))\Big)=0,\\
    &\z=\int_Z f(x,y,u(z))\,dz,\, \zeta^0\geq \int_Z f^0(x,y,u(z))\,dz,\,u(\cd)\in \cU_Z \Big\}.
\end{split}
\end{equation}

Our main result is the following theorem.
\begin{theorem}\label{main}
  Assume $(S1)$ --- $(S5)$, and the following condition hold
  \begin{equation}\label{Cesari}
    \cE(x,y)=\bigcap_{\d>0}\overline{G\cE(x, B_\d(y))},\q \eqae\, (x,y)\in \O\times
\IR,
  \end{equation}
where $\cE(x,y)$ and $G\cE(x,y)$ are defined by \refeq{E} and \refeq{GE},
 $\overline{G\cE(x,y)}$ is the closure of $G\cE(x,y)$ in $\cS^n_+  \times \IR \times \IR$, and $B_\d(y)$ is the ball centered at $y$ with  radius $\d$. Then Problem $(C)$ admits at least one solution.
\end{theorem}

\begin{Remark}
  When the leading term is independent of control variable, i.e. $A(x,u)\equiv A(x)$, Theorem \ref{main}
   is equivalent to the classical existence result of optimal control $($see Theorem~6.4 in Chapter~3 of  \cite{Li}$)$. This fact will follows by Proposition \ref{remark} in Section 3.
\end{Remark}



When dealing  with problems with controls containing in the leading
term, we meet a main difficulty that is to find the state equation
corresponding to the weak limit of state sequence.  This is involved
with the $H$-convergence and $G$-closure problem. It is known that
optimal control usually does not exist for Problem (C) and therefore
to seek optimal relaxed control for Problem (C) is more meaningful
than to seek a solution for Problem (C). Nevertheless, we think this
paper contains some useful ideas for us to get the relaxation of
Problem (C), which will be our forthcoming work. In this paper, we
will give a local representation of G-closure in
 Section 2, which is critical in proving the existence theorem. While Section 3 is devoted to a proof of Theorem
\ref{main} and some propositions.

\section{H-convergence and Local Representation of G-closure }

Now, let us recall the notion of H-convergence. This kind of
convergence was introduced by Murat and Tartar in \cite{Murat}.

\begin{definition}
  A sequence of matrix valued functions $A^\ve(\cd)\in L^\infty(\O;\mcM_{\L,\l})$ is said to H-converge
 to a  matrix valued function $A^*(\cd)\in L(\O;\mcM_{\L,\l})$, if for any right hand side $f\in H^{-1}(\O)$,
 the sequence $y^\ve(\cd)\in H^1_0(\O)$ of weak solutions of
  \beq{eqve1}{\left\{
    \begin{array}{ll}
    \ds -\na\cd(A^\ve(x)\na y^\ve(x))=f,& \eqin \O,\\
    \ds y^\ve(x)=0,& \eqon \pa \O
    \end{array}
    \right.
}
satisfies
\beq{}{\nnb y^\ve(\cd)\wto \by(\cd),\q\wi{\, H_0^1(\O)},}
where $\by(\cd)$ is the weak solution of
\beq{by1}{\left\{
    \begin{array}{ll}
     \ds -\na\cd(A^*(x)\na \by(x))=f,& \eqin \O,\\
    \ds \by(x)=0,& \eqon \pa\O.
    \end{array}
    \right.
}
\end{definition}

The notion of ``H-convergence" makes sense because of the next
compactness proposition (see Theorem 2 in \cite{Murat}).

\begin{prop}\label{compact} 
  For any sequence $A^\ve(\cd)$ of matrices in
  $L^\infty(\O;\mcM_{\L,\l})$,
there exists a subsequence of $A^\ve(\cd)$, $H$-converges to an
$A^*(\cd) \in L^\infty(\O;\mcM_{\L,\l})$.
\end{prop}

This proposition proves the existence of an $H$-limit for a
subsequence of a bounded sequence, but it delivers no explicit
formula for this limit. The next proposition shows that when
$A^\ve(\cd)=A({\cd\over \ve})$ with some periodic matrix valued
function $A(\cd)$,   $A^\ve(\cd)$ $H$-converges to an $H$-limit
defined by an explicit formula (up to solving some corresponding
cell problems). The proposition can be stated as
\begin{prop}\label{periodic}
  Let $A(\cd)\in L_\#^\infty(Z;\mcM_{\L,\l})$. Then
\beq{}{\nnb A\Big({\cd\over \ve}\Big)\Hto A^*\chi_\O(\cd)}
with $A^*\in \mcM_{\L,\l}$ being a constant matrix defined by its
entries
\beq{}{A^*_{ij}=\int_Z A(z)(e_i+\na w_i)\cd e_j \,dz,}
where $\set{w_i}_{1\leq i \leq n}$ is the family of unique solutions
in $H^1_\#(Z)/\IR$ of the cell problems
\beq{}{-\na\cd (A(z)(e_i+\na w_i(z)))=0,\q \eqin \, Z.}
\end{prop}

For a proof of the above proposition, see Theorem 1.3.18 of
\cite{alla2} or Theorem 1.3.1 of \cite{Ben}. The next classical
result (see Theorem 1.3.23 in
\cite{alla2}, for example) shows the fact that 
a general $H$-limit $A^*(\cd)$ can be attained as the limit of a
sequence of periodic homogenized matrices.

\begin{prop}\label{noperiodic}
  Assume  $A^\ve(\cd)\in L^\infty(\O;\mcM_{\L,\l})$   $H$-converge  to
  a limit $A^*(\cd)$. For any $x$ in   $\O$ and
any sufficiently small positive $h>0$, let $A^*_{\ve,h}(\cd)$ be the
periodic homogenized matrix defined by its entries
  \beq{}{\left(A^*_{\ve,h}(x)\right)_{ij}=\int_Z A^\ve(x+hz)(e_i+\na w^i_{\ve,h}(z;x))\cd e_j\,dz,}
  where $\set{w^i_{\ve,h}(\cd;x)}_{1\leq i \leq n}$ is the family of unique solutions in $H^1_\#(Z)/\IR$
  of the cell problems
\beq{}{-\na_z\cd (A^\ve(x+hz)(e_i+\na_z w^i_{\ve,h}(z;x)))=0,\q
\eqin \, Z.}
Then, along a subsequence $h\to 0$,
\beq{E209}{\nnb \lim_{\ve\to 0^+}A^*_{\ve,h}(x)\to A^*(x),\q\eqae\,
x\in
  \O.}
%
\end{prop}

We list some useful properties of $H$-convergence in  follows. For
proofs of these results, see Proposition 1.2.18, Proposition 1.2.22
and Proposition 1.3.44 in \cite{alla2}.

\begin{prop}\label{local}
  Let $A^\ve(\cd)$ and $B^\ve(\cd)$ be two sequences in  $L^\infty(\O;\mcM_{\L,\l})$,
   which $H$-converge to $A^*(\cd)$ and $B^*(\cd)$,
   respectively. Let $\O_0$ be an open subset of $\O$. If
  \beq{}{\nnb A^\ve(x)=B^\ve(x),\q \eqin  \O_0,}
 then
  \beq{}{\nnb A^*(x)=B^*(x),\q \eqin  \O_0.}
\end{prop}

Proposition \ref{local} shows that the value of $H$-limit $A^*(\cd)$
in a region $\O_0$ does not depend on he values of sequence
$A^\ve(\cd)$ outside of this region, which is precisely what we mean
by locality.

\begin{prop}\label{strongH}
  Assume $A^\ve(\cd)\in L^\infty(\O;\mcM_{\L,\l})$
  converge  strongly to  a limit matrix $A^*(\cd)\in L^1(\O; \mcM_{\L,\l})$. Then, $A^\ve(\cd)$   $H$-converges to $A^*(\cd)$ too.

  In particular, if $A^\ve(\cd)\in \mcM_{\L,\l}$ converges to $A^*(\cd)$ almost everywhere in
  $\O$, then $A^\ve(\cd)$  $H$-converges to $A^*(\cd)$.
\end{prop}
This proposition shows that $H$-convergence is weaker than strong
convergence. On the other hand, it is well-known that usually the
weak limit of a sequence $A^\ve(\cd)$ does not equal to its
$H$-limit.

\begin{prop}\label{correct}
  Let \thb{S1} hold. Then there exist  constants $C>0$ and $\d>0$ such that, for any $1\leq p\leq 1+\d$
  and two sequences of  $A^\ve(\cd)$ and $B^\ve(\cd)$ in $L^\infty(\O;\mcM_{\L,\l})$, which H-converge
   to $A^*(\cd)$ and $B^*(\cd)$, respectively, it holds that
  \beq{}{\|A^*(\cd)-B^*(\cd)\|_{L^p(\O)}\leq C\liminf_{\ve\to 0^+}\|A^\ve(\cd)-B^\ve(\cd)\|_{L^p(\O)}.}
\end{prop}

Now define
\beq{}{\cG(A)=\set{P(\cd)\in
L^\infty(\O;\mcM_{\L,\l})\,\Big|~\exists
    u^\ve(\cd)\in\cU,~~\text{s.t.}\, A(\cd,u^\ve(\cd))\Hto A^*(\cd)}.}
We see that $\cG(A)$ is the set of all possible $H$-limits of
$\ds{\set{ A(\cd,u(\cd))}_{u(\cd)\in \cU}}$. A very important
problem called $G$-closure problem is to find out the structure of
$\cG(A)$. Many works devoted to this problem dealt with two-phrase
composite cases(see, for examples, \cite{Cibib}, \cite{Milton} and
\cite{Tartar}). In \cite{Tartar}, a  precise formula of $\cG(A)$ was
given for a special two-phrase case of $A$ taking only $\a I$ and
$\b I$ for some $\b>\a>0$. Unfortunately, in most cases including
usual two-phrase cases, precise knowledge of the $G$-closure are
still lacking.

A local representation of $\cG(A)$ is crucial to our main result. We
give a simple lemma related to Assumption (S3) first.
\begin{lemma}\label{T208} Let $\o(\cd)$ be continuous on $[0,+\infty)$, $\o(0)=0$ and   $F(\cd)\in L^\infty(\O;\IR^m)$.
We have the following results.

\noindent \thb{i}
\begin{equation}\label{E214}
\lim_{h\to 0}\int_Z \o\Big(\big| F(x+hz)-F(x)\big|\Big)\, dz=0,
\qq\eqae x\in \O.
\end{equation}

\noindent \thb{ii} Let $\set{\O_j^k}_{1\leq j\leq k}$ be a family of
measurable decompositions of $\O$ such that:

\thb{a} if $i\ne j$,  $ \ds\O_i^k\bigcap \O_j^k=\emptyset $;

\thb{b} for any $k$, $\ds\bigcup_{1\leq j\leq k} \O_j^k=\O$;

\thb{c}  $\ds \lim_{k\to +\infty}\max_{1\leq j\leq k}\diam (\O_j^k)=
0$.

Then
\begin{equation}\label{E215}
\lim_{k\to \infty}\sum^k_{\ell=1}{1\over |\O^k_\ell |}\int_{\O^k_\ell}
\int_{\O^k_\ell}\o\Big( \big| F(s)- F(x)\big|\Big)\, ds  \, dx=0,
\end{equation}
where $|E|$, $\diam (E)$ denotes the Lebesgue measure and   the
diameter of $E$, respectively.
\end{lemma}
\textbf{Proof.} (i) Let $x\in \O$ be a Lebesgue point of $F(\cd)$
and satisfy $\ds |F(x)|\leq \|F\|_{L^\infty(\O;\IR^m)}$. Then
\begin{equation}\label{E216}\nnb
\lim_{h\to 0}\int_Z  \big| F(x+hz)-F(x)\big| \, dz=0.
\end{equation}
Thus,  as a function of $z$, $\big| F(x+hz)-F(x)\big|$ converges in
measure to $0$  as $h\to 0$. Since
$$
\o\Big(\big| F(x+hz)-F(x)\big|\Big)\leq \max \set{\o (r)\big| 0\leq
r\leq 2\|F\|_{L^\infty(\O;\IR^m)}}  , \qq\eqae z\in Z,
$$
we get \refeq{E214} by Lebesgue's dominated convergence theorem.

(ii) By Remark \ref{R102}, we suppose $\o(\cd)$ is a continuous
module without loss of generality. For any $\ds\F(\cd)\in
C(\bGO;\IR^m)$, we have
\begin{eqnarray}
\nnb & & \sum^k_{\ell=1}{1\over |\O^k_\ell |}\int_{\O^k_\ell}
\int_{\O^k_\ell} \o\Big(\big| F(s)- F(x)\big|\Big)\, ds  \, dx\\
\nnb &\leq & \sum^k_{\ell=1}{1\over |\O^k_\ell |}\int_{\O^k_\ell}
\int_{\O^k_\ell} \Big\{\o\Big(\big| F(s)-
\F(s)\big|\Big)+\o\Big(\big| \F(s)- \F(x)\big|\Big)+\o\Big(\big|
\F(x)- F(x)\big|\Big)\Big\}\, ds \, dx\\
\nnb &=&  2\int_\O  \o\Big(\big| F(x)- \F(x)\big|\Big) \, dx
+\sum^k_{\ell=1}{1\over |\O^k_\ell |}\int_{\O^k_\ell}
\int_{\O^k_\ell}  \o\Big(\big| \F(s)- \F(x)\big|\Big) \, ds \, dx.
\end{eqnarray}
Consequently, if we set $F(x)=0$ for $x\not\in \O$ and choose
$$
\F(x)=\int_Z F(x+hz)\, dz,
$$
we have $\F\in C(\bGO;\IR^m)$. Consequently, it follows easily from
the uniform continuity of $\F$ and the assumption (c) that
$$
\lim_{k\to +\infty}\sum^k_{\ell=1}{1\over |\O^k_\ell
|}\int_{\O^k_\ell} \int_{\O^k_\ell}  \o\Big(\big| \F(s)-
\F(x)\big|\Big) \, ds \, dx=0.
$$
Thus,
\begin{eqnarray}\label{E218}
\nnb & & \limsup_{k\to +\infty}\sum^k_{\ell=1}{1\over |\O^k_\ell
|}\int_{\O^k_\ell}
\int_{\O^k_\ell} \o\Big(\big| F(s)- F(x)\big|\Big)\, ds  \, dx\\
\nnb &\leq &   2\int_\O  \o\Big(\big| F(x)- \int_Z F(x+hz)\,
dz\big|\Big) \, dx\\
\nnb  &\leq &  2\int_\O  \o\Big(\int_Z \big|F(x+hz)-F(x)\big|\, dz
\Big) \, dx, \qq\all h>0.
\end{eqnarray}
Therefore,
\begin{eqnarray}\label{E219}
\nnb & & \limsup_{k\to +\infty}\sum^k_{\ell=1}{1\over |\O^k_\ell
|}\int_{\O^k_\ell}
\int_{\O^k_\ell} \o\Big(\big| F(s)- F(x)\big|\Big)\, ds  \, dx\\
\nnb &\leq &     2\lim_{h\to 0}\int_\O  \o\Big(\int_Z
\big|F(x+hz)-F(x)\big|\, dz \Big) \, dx=0.
\end{eqnarray}
We get the proof. \endpf

Now, we will give a local representation of $\cG(A)$.
\begin{theorem}\label{Gclosure}
  Assume \thb{S1}---\thb{S3} hold. Then the $G$-closure set $\cG(A)$ is characterized by
  \beq{}{\cG(A)=\set{P(\cd)\in L^\infty(\O;\mcM_{\L,\l})~\Big|~~P(x)\in \ol{G_x(A)},\,\eqae{\,x\in\O}} }
with $G_x(A)$ being defined by
\beq{Gx}{\begin{split}
     G_x(A)=\Big\{Q\in \cS^n_+~\Big|&~\exists u\in\cU_Z,~~\text{s.t.}~~Q_{ij}=\int_Z A(x,u(z))(e_i+\na w^i (z;x))\cd e_j\,dz,
    \\ & \text{where}~~w^i(\cd;x)\in H^1_\#(Z), ~~ \na_z\cd\Big(A(x,u(z))(e_i+\na w^i (z;x))\Big)=0 \Big\}.
\end{split}
}
\end{theorem}
\begin{Remark} It is easy to see that in \refeq{Gx}, $G_x(A)$ can be rewritten as

\begin{equation}
\begin{split}
     G_x(A)=\Big\{Q\in \cS^n_+~\Big|&~\exists u\in\cU_Z,~~\text{s.t.}~~Q_{ij}=\int_Z A(x,u(z))(e_i+\na w^i (z;x))\cd (e_j+\na w^j(z))\,dz,
    \\&~~\text{where}~~w^i(\cd,;x)\in H^1_\#(Z), ~~ \na_z\cd\Big(A(x,u(z))(e_i+\na w^i (z;x))\Big)=0
\Big\},
\end{split}
\end{equation}
which implies  $\ds G_x(A)\subseteq \mcM_{\L,\l}$.
\end{Remark}
\textbf{Proof of Theorem 2.9.} Denote
 \beq{}{\nnb \cP(A)=\set{P(\cd)\in L^\infty(\O;\mcM_{\L,\l})~\Big|~~P(x)\in \ol{G_x(A)},\,\eqae{\,x\in\O}}. }

We need to show $\cG(A)=\cP(A)$.

We prove $\cG(A)\subseteq\cP(A)$ first. Assume $A^*(\cd)\in \cG(A)$.
Then there exists a sequence $u^\ve(\cd)\in \cU$, such that as
$\ve\to 0^+$,
\beq{}{\nnb A(\cd,u^\ve(\cd))\Hto A^*(\cd).}
By Proposition \ref{noperiodic}, along a subsequence $h\to 0$,
\beq{E217}{A^*_{ij}(x)=\lim_{h\to 0}\lim_{\ve\to
    0}(\tiA^*_{h,\ve}(x))_{ij}\qq\eqae x\in \O,}
where $\tiA^*_{h,\ve}(\cd)$ is defined by
 \beq{tiA}{(\tiA^*_{h,\ve}(x))_{ij}=\int_Z
 A(x+hz,u^\ve(x+hz))(e_i+\na_z
    \tiw^i_{h,\ve}(z;x))\cd e_j\,dz }
with $\tiw^i_{h,\ve}(\cd;x)\in H^1_\#(Z)/\IR$ being the unique
$Z$-periodic solution of
\beq{tiw}{\na_z\cd\Big(A(x+hz,u^\ve(x+hz))(e_i+\na_z
\tiw^i_{h,\ve}(z;x)) \Big)=0.}
On the other hand, define $A^*_{h,\ve}(\cd)$ by 
\beq{A}{(A_{h,\ve}^*(x))_{ij}=\int_Z A(x,u^\ve(x+hz))(e_i+\na_z
w^i_{h,\ve}(z;x))\cd e_j\,dz }
with $w^i_{h,\ve}(\cd;x)\in H^1_\#(Z)/\IR$ being the unique
$Z$-periodic solution of
\beq{w}{\na_z\cd\Big(A(x,u^\ve(x+hz))(e_i+\na_z w^i_{h,\ve}(z;x))
\Big)=0.}
Then, combining \refeq{tiw} with \refeq{w}, we get
\begin{eqnarray}\label{E222}
    \nnb&&\na_z\cd\Big(A\big(x,u^\ve(x+hz)\big)\big(\na \tiw^i_{h,\ve}(z;x)-\na w^i_{h,\ve}(z;x)\big)\Big)\\
    &=&\na_z\cd\Big[\Big(A\big(x,u^\ve(x+hz)\big)-A\big(x+hz,u^\ve(x+hz)\big)\Big)\big(e_i+\na_z \tiw^i_{h,\ve}(z;x)\big)\Big].
\end{eqnarray}
Multiplying \refeq{E222} by $\ds
\tiw^i_{h,\ve}(z;x)-w^i_{h,\ve}(z;x) $ and using integration by
part, we get from the periodicities of $\tiw^i_{h,\ve}(\cd;x)$ and
$w^i_{h,\ve}(\cd;x)$ that
\begin{eqnarray}
    \nnb &&\int_Z A\big(x,u^\ve(x+hz)\big)\big(\na \tiw^i_{h,\ve}(z;x)-\na w^i_{h,\ve}(z;x)\big)\cd
    \big(\na \tiw^i_{h,\ve}(z;x)-\na w^i_{h,\ve}(z;x)\big)\,dz\\
    \nnb &=&\int_Z \Big(A(x,u^\ve(x+hz))-A(x+hy,u^\ve(x+hz))\Big)\\
    && \q\big(e_i+\na_z \tiw^i_{h,\ve}(z;x)\big)
    \cd\big(\na_z \tiw^i_{h,\ve}(z;x)-\na_z w^i_{h,\ve}(z;x)\big)\,dz.
\end{eqnarray}
Then the ellipticity of $A$ yields
\begin{eqnarray}
    \nnb &&\l\|\na w^i_{h,\ve}(\cd;x)-\na \tiw^i_{h,\ve}(\cd;x)\|_{L^2(Z)}\\
    \nnb &\leq&\left\{\int_Z\left|\Big(A(x,u^\ve(x+hz))-A(x+hz,u^\ve(x+hz))\Big)\big(e_i+\na w^i_{h,\ve}(z;x)\big)\right|^2\,dz\right\}^{1/2}.
\end{eqnarray}
By \refeq{w} and Meyers' theorem (see \cite{Meyers}, see also
Theorem 1.3.41 and Remark 1.3.42  in \cite{alla2}), there exist
constants $p>2$ and $C>0$, both dependent only on $\l$, $\L$ and
$\O$, such that
\begin{equation}\label{E225}
\| \na \tiw^i_{h,\ve}(\cd;x) \|_{L^p(Z)}\leq C, \q \|\na
w^i_{h,\ve}(\cd;x) \|_{L^p(Z)}\leq C.
\end{equation}
 Thus
 \beq{wtiw}{\|\na w^i_{h,\ve}(\cd;x)-\na
    \tiw^i_{h,\ve}(\cd;x)\|_{L^2(Z)}\leq
    C\Big(\int_Z|A(x,u^\ve(x+hy))-A(x+hy,u^\ve(x+hy))|^q\,dy\Big)^{1/q},}
where $\ds{{1\over p}+{1\over q}={1\over 2}}$. Then it follows from
\refeq{wtiw} and (S3) that
\begin{eqnarray}\label{use}
       \nnb & & |(\tiA_{h,\ve}^*(x))_{ij}-(A_{h,\ve}^*(x))_{ij}|\\
        \nnb &\leq & \int_Z|A(x,u^\ve(x+hz))(\na_z \tiw^i_{h,\ve}(z;x)-\na_z w^i_{h,\ve}(z;x))|\,dz\\
        \nnb && +\int_Z\left|\Big(A(x+hz,u^\ve(x+hz))-A(x,u^\ve(x+hz))\Big)(e_i+\na_z \tiw^i_{h,\ve}(z;x))\right|\,dz\\
        \nnb &\leq & C\|\na  \tiw^i_{h,\ve}(\cd;x)-\na w^i_{h,\ve}(\cd;x))\|_{L^2(Z)} \\
        \nnb & & +         C\Big(\int_Z |A(x+hz,u^\ve(x+hz))-A(x,u^\ve(x+hz))|^q\,dz\Big)^{1/q}
        \|e_i+\na  \tiw^i_{h,\ve}(\cd;x) \|_{L^p(Z)}\\
       \nnb  &\leq& C \Big(\int_Z|A(x+hz,u^\ve(x+hz))-A(x,u^\ve(x+hz))|^q\,dz\Big)^{1/q}\\
      &\leq&C \left\{\int_Z \Big[\o\Big(\big|F(x+hz)-F(x)\big|\Big)\Big]^q\,dz\right\}^{1/q}.
\end{eqnarray}
Thus, by Lemma \ref{T208}, we get from \refeq{E217} and \refeq{use}
that along a subsequence $h\to 0$,
\beq{uvex}{A^*_{ij}(x)=\lim_{h\to 0}\lim_{\ve\to
    0}(A^*_{ h,\ve}(x))_{ij}, \qq\eqae x\in \O.}
Noting that $A^*_{h,\ve}(x)\in G_x(A)$, we get $A^*(x)\in
\ol{G_x(A)}$, $\eqae x\in\O$. Therefore $\cG(A)\subseteq\cP(A)$.

Next, we turn to prove $\cP(A)\subseteq\cG(A)$.  Let
$\set{\O_j^k}_{1\leq j\leq k}$ be a family of measurable
decompositions of $\O$ satisfying (a)--(c) in Lemma \ref{T208}.
Denote by $\chi_j^k(\cd)$ the characteristic function of $\O_j^k$.

We will  show the result in three steps.

\textbf{Step I.} Assume  $A(x,u)\equiv A(u)$.

Denote $G(A)\equiv G_x(A)$ since $G_x(A)$ is independent of $x$ in
this case.

For any $\ds \tiA   \in G(A)$, we have $u(\cd)\in \cU_Z$ such that
$$
 \tiA_{ij}=\int_Z A(u(z))(e_i+\na w^i (z))\cd e_j\,dz,
$$
where $w^i(\cd)\in H^1_\#(Z)/\IR$ solves
$$
\na \cd\Big(A(u(z))(e_i+\na w^i (z))\Big)=0.
$$
By Proposition \ref{periodic}, $\ds A\big(u\big({\cd\over
\ve}\big)\big)\Hto \tiA\chi_\O(\cd)$ as $\ve\to 0^+$. Thus, $\tiA
\chi_\O(\cd)\in  \cG(A)$. We denote this result simply by $\ds
G(A)\subseteq \cG(A)$. Obviously, we can get $\ds \ol{G(A)}\subseteq
\cG(A)$ immediately.

Let $A^*(\cd)\in \cP(A)$. Then
\beq{}{\nnb A^*(x)\in \overline{G(A)},\q\eqae  x\in\O. }
Define
\beq{}{\nnb \hA^k_j={1\over|\o_j^k|}\int_{\o^k_j}A^*(x)\,dx,\q
\hA_k(\cd)=\sum_{j=1}^k \hA^k_j\chi_j^k(\cd).}
 Then
\beq{hA}{\hA^k(\cd)\to A^*(\cd)\qq\si{\,L^p(\O)}, \qq \all 1\leq
p<\infty.}
Denote
 \beq{projection}{
\cA^k_j\equiv\set{Q\in \ol{G(A)}\Big|\,
    \big|Q-\hA^k_j\big|=\inf_{P\in \ol{G(A)}}\big|P-\hA^k_j\big|},\q 1\leq j \leq k;\, k=1,2,3,\ldots.}
Since $\overline{G(A)}$ is closed, $\cA^k_j$ is always nonempty.
Thus, we can select a constant matrix $A^k_j$ from $\cA^k_j$. Define
\beq{}{A_k(\cd)=\sum_{j=1}^n A^k_j\chi_j^k(\cd).}
Since for almost all $x\in \O^k_j$, $A^*(x)\in \ol{G(A)}$, by
\refeq{projection}, there is
\begin{eqnarray*}
      &   &|A_k(x)-A^*(x)|\\
    &\leq &|A^k_j-\hA^k_j|+|\hA^k_j-A^*(x)|\\
    &\leq &2|\hA^k_j-A^*(x)|.
\end{eqnarray*}
Thus, by \refeq{hA},
\beq{}{A_k(\cd)\to A^*(\cd),\qq\si{\,L^1(\O)}.}
Consequently, by  Proposition \ref{strongH}, we have
\beq{H}{A_k(\cd)\Hto A^*(\cd).}
The advantage of replacing $\ds \hA_k$ by  $\ds A_k$ is that
 we have
 $ A^k_j \in \ol{G(A)}\subseteq  \cG(A)$ while we do not always have  $ \hA^k_j \in
 \ol{G(A)}$.

Then, by Proposition \ref{local} (local property),
$A_k(\cd)\in\cG(A)$. Finally, by \refeq{H}, $A^*(\cd)\in \cG(A)$.
That is, $\cP(A)\subseteq \cG(A)$.

\textbf{Step II.}  Assume $\ds  A(x,u)=\sum_{1\leq j\leq k}
A_j(u)\chi^k_j(x) $.        

 By what we have proved in Step I and the local property of $H$-convergence, we can see that  $\cP(A)\subseteq
\cG(A)$ holds in this case.

\textbf{Step III.} General cases. Let $A^*(\cd)\in \cP(A)$. Then
$A^*(x)\in \ol{G_x(A)}$, $\eqae x\in \O$.  We want to prove
$A^*(x)\in \cG(A)$. Without loss of generality, we can suppose that
\begin{equation}\label{E238}
A^*(x)\in G_x(A), \qq \all x\in \O.
\end{equation}
Define
\beq{}{\nnb\left\{\begin{array}{l} \ds A_k(x,u)=\sum_{j=1}^k
A^k_j(u)\chi_j^k(x),\\
\ds
   A^k_j(u)={1\over|\O_j^k|}\int_{\O^k_j}A(s,u)\,ds,\end{array}\right. \qq (x,u)\in \O\times
   U.}
While $A_k^*(\cd)$ is  a measurable selection of the projection of
$A^*(x)$ on $\overline{G_x(A_k)}$, i.e., $A_k^*(\cd)$ is measurable
and
$$
\big|A_k^*(x)-A^*(x)\big|=\inf_{P\in
\ol{G_x(A_k)}}\big|P-A^*(x)\big|, \q A_k^*(x)\in \ol{G_x(A_k)},
$$
where $G_x(A_k)$ is defined by \refeq{Gx}. By Filippov's lemma (see
\cite{Filippov},  or Corollary 2.26 of Chapter 3 in \cite{Li}), such
an $A_k^*(\cd)$ exists.

Now we will show that
\beq{conv1}{A^*_k(\cd)\to A^*(\cd),\qq\si\, L^1(\O;\cS^n_+).}

By \refeq{E238}, there exists a $u(\cd)\in\cU_Z$, such that
\beq{}{\nnb(A^*(x))_{ij}=\int_Z A(x,u(z))(e_i+\na_z w^i(z;x))\cd
e_j\,dz,}
where for any $x\in \O$, $w^i(\cd;x)\in H^1_\#(Z) /\IR$ is the
solution of
\beq{}{\nnb \na_z\cd\Big(A(x,u(z))(e_i+\na_z w^i(z;x))\Big)=0.}
Next, we can define $\tiA^*_k(\cd)$ by
\beq{}{\nnb(\tiA^*_k(x))_{ij}=\int_Z A_k(x,u(z))(e_i+\na_z
w_{k}^i(z;x))\cd e_j\,dz,}
where $w_{k}^i(\cd;x)\in H^1_\#(Z)/\IR$ is the solution of
\beq{}{\nnb\na_z \cd\Big(A_k(x,u(z))(e_i+\na_z w_{k}^i(z;x))\Big)=0.}
This means that $\tiA^*_k(x)\in\ol{G_x(A_k)}$. Thus, similar to the
proof of \eqref{use}, we have
\begin{eqnarray}\label{E245}
     \nnb   && |(A^*_k(x))_{ij}-(A^*(x))_{ij}|\\
\nnb    &\leq &|(\tiA^*_k(x))_{ij}-(A^*(x))_{ij}|\\
\nnb &\leq & \int_Z |A(x,u(z))(\na_z w^i(z;x)-\na_z w_k^i(z;x))|\,dz \\
\nnb & & +\int_Z|\big(A_k (x,u(z))-A(x,u(z))\big)(e_i+\na_z w_k^i(z;x))|\,dz \\
  &  \leq & \Big(\int_Z |A_k(x,u(z))-A(x,u(z))|^q\,dz\Big)^{1/q}.
\end{eqnarray}
By Lebesgue's dominated convergence theorem, we deduce
\beq{E246}{\lim_{k\to
\infty}|(A^*_k(x))_{ij}-(A^*(x))_{ij}|=0,\qq\eqae{x\in\O},}
which proves \refeq{conv1}.

Furthermore, noting that $A_k^*(\cd)$ is piecewise constant  and
$A_k^*(\cd)\in \cP(A_k)$, by Step II, $A_k^*(\cd)\in \cG(A_k)$. Then
there exists  $u^{k,j}(\cd)\in \cU$, such that
\beq{}{A_k(\cd,u^{k,j}(\cd))\Hto  A^*_k(\cd), \qq (j\to +\infty).}
By Proposition \ref{compact}, we can suppose that
\beq{}{A(\cd,u^{k,j}(\cd))\Hto  A_k(\cd), \qq (j\to +\infty).}
By Proposition \ref{correct}, we obtain
\begin{eqnarray}\label{conv2}
\nnb &  & \|A^*_k(\cd)-A_k(\cd)\|_{L^1(\O)}\\
\nnb  &\leq & C\liminf_{j\to +\infty}\|A_k(\cd,u^{k,j}(\cd))-A(\cd,u^{k,j}(\cd))\|_{L^1(\O)}\\
\nnb  &=& C \liminf_{j\to +\infty}\int_\O\Big|
\sum^k_{\ell=1}{1\over |\O^k_\ell |}
\int_{\O^k_l}\Big(A(s,u^{k,j}(x))-A(x,u^{k,j}(x))\Big)\, ds \, \chi^k_\ell(x)\Big|\, dx\\
\nnb  &=& C \liminf_{j\to +\infty}\sum^k_{\ell=1} {1\over |\O^k_\ell
|}\int_{\O^k_\ell} \Big|
\int_{\O^k_\ell}\Big(A(s,u^{k,j}(x))-A(x,u^{k,j}(x))\Big)\, ds  \Big|\, dx\\
\nnb&\leq & C \sum^k_{\ell=1}{1\over |\O^k_\ell |}\int_{\O^k_\ell}
\int_{\O^k_\ell} \o\Big(\big|F(s)-F(x)\big|\Big)\, ds  \, dx .
\end{eqnarray}
Thus it follows from Lemma \ref{T208} that,
$$
\lim_{k\to +\infty} \|A^*_k(\cd)-A_k(\cd)\|_{L^1(\O)}=0.
$$
Combining the the above with \refeq{E246}, we get
$$
\lim_{k\to +\infty} \|A_k(\cd)-A^*(\cd)\|_{L^1(\O)}=0.
$$
Consequently,
\beq{}{A_k(\cd)\Hto A^*(\cd).}
It follows from $A_k(\cd)\in \cG(A)$  that   $A^*(\cd)\in \cG(A)$.
This ends the proof.
\endpf

\section{Proof of the Main Theorem}
In this section, we will prove our main result. Before that, we need
to show three lemmas. The first is  about the well-posedness and
regularity of state equation \refeq{orgin}.

\begin{lemma}\label{wellpose}
  Let \thb{S1}---\thb{S4} hold. Then for any $u(\cd)\in \cU$, \refeq{orgin} admits a
  unique weak solution $y(\cd)\in H^1_0(\O)\cap L^\infty(\O)$. Furthermore, there exists a constant $R>0$,
  independent of $u(\cd)$, such that
  \beq{}{\|y\|_{H^1_0(\O)}+\|y\|_{L^\infty(\O)}\leq R.}
\end{lemma}
The existence of a weak solution to \refeq{orgin} in $H^1_0(\O)$
together with the $H^1_0(\O)$-norm  estimate follows easily from the
variational structure of \refeq{orgin}, while the uniqueness of the
weak solution follows from (S3) and \refeq{fy}. The boundedness of
weak solution in $L^\infty(\O)$ follows from standard De Giorgi
iteration.

In order to proof our main theorem, we need another lemma.
\begin{lemma}\label{lemma}
  Assume $A^\ve(\cd)\in L^\infty(\O,\mcM(\L,\l))$ and
\beq{}{\nnb A^\ve(\cd)\Hto  A^*(\cd),\qq (\ve\to 0^+).}
Moreover $f^\ve(\cd), f(\cd)\in L^2(\O)$ and  $f^\ve(\cd)\wto
f(\cd)$ weakly in $L^2(\O)$. Let
  $y^\ve(\cd)\in H_0^1(\O)$ be the weak solution of
\beq{eqve}{\left\{
    \begin{array}{ll}
    \ds -\na\cd(A^\ve(x)\na y^\ve(x))=f^\ve(x),& \eqin\, \O,\\
     y^\ve(x)=0,& \eqon \pa \O.
    \end{array}
    \right.
}
Then
\beq{}{\nnb y^\ve(\cd)\wto \by(\cd)\qq\wi{H_0^1(\O)},}
where $\by(\cd)$ is the weak solution of
\beq{by}{\left\{
    \begin{array}{ll}
   \ds -\na\cd(A^*(x)\na \by(x))=f(x),& \eqin{\,\O},\\
    \ds \by(x)=0,& \eqon{\q\pa\O}.
    \end{array}
    \right.
}
\end{lemma}
\textbf{Proof.}  Set $h^\ve(\cd)=f^\ve(\cd)-f(\cd)$. Then $\ds
\|h^\ve(\cd)\|_{L^2(\O)}$ is bounded. Let
 $z^\ve(\cd)\in H^1_0(\O)$ be the weak solution of
\beq{}{\left\{
    \begin{array}{ll}
    \ds -\na\cd(A^\ve(x)\na z^\ve(x))=h^\ve(x),& \eqin{\,\O},\\
    \ds z^\ve(x)=0,& \eqon{\,\pa\O}.
    \end{array}
    \right.
}
We have
\begin{equation}\nnb
\begin{split}
  & \l\int_\O|\na z^\ve(x)|^2\,dx\leq\int_\O A^\ve(x)\na z^\ve(x)\cd\na z^\ve(x)\,dx\\
  =&\int_\O z^\ve(x)h^\ve(x)\,dx
  \leq C\|z^\ve(\cd)\|_{L^2(\O)}.
\end{split}
\end{equation}
Hence $\ds \|z^\ve(\cd)\|_{H^1_0(\O)}$ is bounded. Then along a
subsequence $\ve\to 0^+$,
$$
z^\ve(\cd)\to\bz(\cd),\qq\wi{\,H^1_0(\O)},\,  \si{L^2(\O)}.
$$
Consequently, along a subsequence $\ve\to 0^+$,
\begin{equation}\nnb
\begin{split}
  & \l\int_\O|\na z^\ve(x)|^2\,dx\leq\int_\O z^\ve(x)h^\ve(x)\,dx\to
  0,
\end{split}
\end{equation}
which means
$$
z^\ve(\cd)\to\bz(\cd)=0,\qq\ si{\,H^1_0(\O)}.
$$
 Moreover, we can get that $z^\ve(\cd)$
itself converges to $0$ strongly in $H^1_0(\O)$.

Since
\beq{}{\left\{
    \begin{array}{ll}
    \ds -\na\cd \Big(A^\ve(x)\na \big(y^\ve(x)-z^\ve(x)\big)\Big)=f(x),& \eqin \O,\\
    \ds y^\ve(x)-z^\ve(x)=0,& \eqon \pa \O,
    \end{array}
    \right.}
\beq{}{\nnb y^\ve(\cd)-z^\ve(\cd)\wto \by(\cd),\q \wi{\,H_0^1(\O)}.}
Thus
$$
y^\ve(\cd)=(y^\ve(\cd)-z^\ve(\cd))+z^\ve(\cd)\wto
\by(\cd),\qq\wi{\,H_0^1(\O)}.
$$
This ends the proof.
\endpf

The third lemma is about relaxed control defined by finite-additive
probability measures. Denote $C(U)$ the bounded continuous function
space on $U$, and $\cM(U)$ the space of all regular bounded finitely
additive measures on $U$. Moreover, denote
   \begin{equation*}
     \cM^1_+(U)=\set{\m\in \cM(U) \Big| \m \text{\, is nonnegative and}\, \m(U)=1}
   \end{equation*}
   and
   \begin{equation*}
     \cR(\O, U)=\set{\s:\O\to \cM^1_+(U)\Bigg|\,x\mapsto\int_Uh(v)\s(x)(dv) \text{\, is measurable}, \all h\in C(U)}.
   \end{equation*}
  Let $C(U)^*$ and $L^1(\O;C(U))^*$ be the dual spaces of $C(U)$ and $L^1(\O;C(U))$, respectively. We regard $\cM^1_+(U)$ and $\cR(\O, U)$ as subspace of $C(U)^*$ and $L^1(\O;C(U))^*$ by setting
   \begin{equation}\label{CU}
     \m(h)\triangleq\int_U h(v)\m(dv),\qq \all h\in C(U),
   \end{equation}
   and
   \begin{equation}\label{LCU}
     \s(g)\triangleq\int_\O\int_U h(x,v)\s(x)(dv),\qq \all g\in L^1(\O; C(U)).
   \end{equation}
   By Theorems 12.2.11 and 12.4.6 in \cite{Fatto}, \refeq{CU} and \refeq{LCU} are well defined. Thus we denote $\s_k(\cd)\Rto \s(\cd)$  if
   \begin{equation}
      \lim_{k\to\infty}\int_\O\int_U h(x,v)\s_k(x)(dv)\,dx=\int_\O\int_U h(x,v)\s(x)(dv)\,dx,\qq \all h\in L^1(\O;C(U)).
   \end{equation}
We have (see Theorem 12.5.9 in \cite{Fatto}):
\begin{lemma}\label{lemma303}
  Assume $(S1)$--- $(S2)$ hold. Let $u_k(\cd)$ be a sequence in
  $\cU$. Then there is a subsequence   of $u_k(\cd)$,
  still denote by itself, such that
$$
\d_{u_k(\cd)}\xrightarrow{\cR}\s(\cd)
$$
for some $\s(\cd)\in \cR(\O, U)$, i.e.
\begin{equation}\label{sconv1}
  \lim_{k\to\infty}\int_\O h(x,u_k(x))\,dx=\int_\O\int_U h(x,v)\s(x)(dv)\,dx,\qq \all h\in L^1(\O;C(U)).
\end{equation}
\end{lemma}

Now we are at the position to prove Theorem \ref{main}.

\textbf{Proof of Theorem 1.1}. Let $ u_k(\cd)\in \cU$ be a
minimizing sequence of Problem (C),  $ y_k(\cd) $ be the
corresponding state sequence. Then
 \beq{}{\nnb\|y_k(\cd)\|_{H^1_0(\O)}+\|y_k(\cd)\|_{L^\infty(\O)}\leq R.}
Thus, along a subsequence,
\beq{yconv}{y_k(\cd)\to\by(\cd)\qq\wi{\, H^1_0(\O)}, \q \eqae \eqin
\O}
for some $\by(\cd)\in H^1_0(\O)\cap L^\infty(\O)$.

By Proposition \ref{compact}, there exists an $A^*(\cd)\in
L^\infty(\O; \mcM(\L,\l))$ and a subsequence of $u_k(\cd)$, still
denoted by $ u_k(\cd) $, such that
\beq{}{\nnb A( \cd,u_k(\cd))\Hto  A^*(\cd).}
Then by \refeq{uvex}, along a subsequence $h\to 0$,
\begin{equation}\label{Aae}
  A^*_{ij}(x)=\lim_{h\to 0}\lim_{k\to \infty}\int_Z A(x,u_k(x+hz))(e_i+w_{h,k}^i(z;x))\cd e_j\,dz.
\end{equation}
where $w_{h,k}^i(\cd;x)\in H^1_\#(Z)$ is the Z-periodic solution of
\begin{equation*}
  \na\cd\Big(A(x,u_k(x+hz))(e_i+w_{h,k}^i(z;x))\Big)=0.
\end{equation*}

On the other hand, since
\beq{}{\nnb |f(x,y_k(x),u_k(x))|\leq M_R,}
we can suppose that
\beq{}{f(\cd,y_k(\cd),u_k(\cd))\wto \barf(\cd), \qq\wi{\, L^2(\O)} }
for some $\barf(\cd)\in L^\infty (\O)$.

   In order to characterize $\barf$  precisely, it is useful to use relax controls defined by finite-additive
   measures. By Lemma \ref{lemma303}, we can suppose that
$$
\d_{u_k(\cd)}\xrightarrow{\cR}\s(\cd)
$$
for some $\s(\cd)$ in $\cR(\O, U)$. That is,
\begin{equation}\label{sconv}
  \lim_{k\to\infty}\int_\O h(x,u_k(x))\,dx=\int_\O\int_U h(x,v)\s(x)(dv)\,dx,\qq \all h\in L^1(\O;C(U)).
\end{equation}
In particular, for any $g\in L^2(\O)$,
$$
  \lim_{k\to\infty}\int_\O f(x,\by(x),u_k(x))\,dx=\int_\O\int_U f(x,\by(x),v)\s(x)(dv)\,dx.
$$
That is,
$$
  f(x,\by(x),u_k(x))\wto \int_U
  f(x,\by(x),v)\s(x)(dv),\qq\wi{L^2(\O)}.
$$
On the other hand, by (S4) and \refeq{yconv},
\begin{eqnarray}
    \nnb&     &\left|  f(x,y_k(x),u_k(x))-  f(x,\by(x),u_k(x))\right|\\
    \nnb&\leq &  M_R  \left|y_k(x)-\by(x)\right|   \to 0,\q
    (k\to\infty).
\end{eqnarray}
Therefore,
\begin{equation*}
  f(x,y_k(x),u_k(x))\wto \int_U f(x,\by(x),v)\s(x)(dv), \q\wi \, L^2(\O),
\end{equation*}
i.e.
\begin{equation}\label{barf}
   \barf(x)=\int_U f(x,\by(x),v)\s(x)(dv).
\end{equation}
Furthermore, define $u_k(x)=0$ if $x\not\in\O$, then for almost all
$x \in \O$,
\begin{equation}\label{f1}
\begin{split}
     & \left|\int_Z f(x,y_k(x),u_k(x+hz))\,dz-\int_Z f(x,\by(x),u_k(x+hz))\,dz\right|\\
\leq & M_R|y_k(x)-\by(x)|\to 0,\qq\,(k\to \infty)
\end{split}
\end{equation}
and
\begin{eqnarray}\label{f2}
\nnb &  &\lim_{k\to \infty}\int_Z f(x,\by(x),u_k(x+hz))\,dz\\
\nnb &= &\lim_{k\to \infty} {1\over h^n}\int_\O f(x,\by(x), u_k(z))\chi_{x+hZ}(z)\,dz\\
\nnb &= &{1\over h^n}\int_\O \int_U f(x,\by(x), v)\chi_{x+hZ}(z)\s(z)(dv)\,dz\\
&= &\int_Z \int_U f(x,\by(x), v)\s(x+hz)(dv)\,dz.
\end{eqnarray}
Combing \refeq{barf}, \refeq{f1} and \refeq{f2}, we obtain
\begin{equation}\label{fae}
  \lim_{h\to 0}\lim_{k\to \infty}\int_Zf(x,y_k(x),u_k(x+hz))\,dz=\barf(x),\q\eqae\, x\in \O.
\end{equation}
In addition, we define
\begin{equation*}
  f^0_{k,h}(x)=\int_Z f^0(x, y_k(x), u_k(x+hz))\,dz,
\end{equation*}
\begin{equation}\label{f0ae}
  \barf^0(x)=\liminf_{h\to 0}\liminf_{k\to \infty}f^0_{k,h}(x),\q \eqae \, x\in \O.
\end{equation}
Then, combing \refeq{Aae}, \refeq{fae} with \refeq{f0ae}, we obtain
that along a subsequence $h\to 0$,
\begin{equation}
  \left(
    \begin{array}{c}
      A^*_{ij}(x) \\
      \barf(x) \\
      \barf^0(x) \\
    \end{array}
  \right)
  =\lim_{h\to 0}\lim_{k\to \infty}
  \left(
  \begin{array}{c}
   \ds \int_Z A(x,u_k(x+hz))(e_i+w_{h,k}^i(z))\cd e_j\,dz\\
   \ds \int_Zf(x,y_k(x),u_k(x+hz))\,dz \\
   \ds \int_Zf^0(x,y_k(x),u_k(x+hz))\,dz\\
   \end{array}
   \right).
\end{equation}
Thus,
\beq{eq11}{\left(A^*(x),\barf(x),\barf^0(x)\right)\in \bigcap_{\d>0}
\overline{G\cE(x,B_\d(\by(x)))},\q \eqae{\,x\in \O}. }
By \refeq{Cesari},
\beq{eq1}{\left(A^*(x),\barf(x),\barf^0(x)\right)\in
 \cE(x,\by(x)),\q \eqae{\,x\in \O}. }
Define
\beq{}{\nnb
g(x,u)=|A(x,u)-A^*(x)|+|f(x,\by(x),u)-\barf(x)|+[f^0(x,\by(x),u)-\barf^0(x)]^+,}
where $a^+$ denote the positive part of a real number $a$. Then,
$g(x,u)$ is measurable in $x$ and continuous in $u$. It follows from
\refeq{eq1} and the definition of $\cE(x,\by(x))$ 
that $0\in g(x, U)$.
 By Filippov's lemma, there exists a $\bu(\cd)\in \cU$, such that
\beq{}{\nnb \left\{\begin{array}{l}
    A^*(x)=A(x,\bu(x)),\\
    \barf(x)=f(x,\by(x),\bu(x)),\\
    \barf^0(x)\geq f^0(x,\by(x),\bu(x)),
    \end{array}\qq\eqae x\in \O.
    \right.}
Consequently,  $\by(\cd)$ is the weak solution of
\beq{}{\nnb \left\{
    \begin{array}{ll}
    \ds -\na\cd(A(x,\bu(x))\na \by(x))=f(x,\by(x),\bu(x)),& \eqin \O,\\
    \ds \by(x)=0,& \eqon \pa\O.
    \end{array}
    \right.
    }
Finally, by Fatou's lemma,
\begin{eqnarray*}
  & &  J(\bu(\cd))= \int_\O f^0(x,\by(x),\bu(x))\,dx\\
  &\leq &\int_\O \barf^0(x)\,dx
    =   \int_\O\liminf_{h\to 0}\liminf_{k\to 0}f^0_k(x)\,dx\\
  &\leq &\liminf_{h\to 0}\liminf_{k\to 0}\int_\O f^0_k(x)\,dx
  \leq  \liminf_{h\to 0}\liminf_{k\to 0}\int_\O \int_Z f^0(x,y_k(x),u_k(x+hz))\,dz \,dx\\
  &=&\lim_{k\to 0} J(u_k(\cd))
  = \inf_{u(\cd)\in\cU}J(u(\cd)).
\end{eqnarray*}
This means that $\bu(\cd)$ is a solution of Problem (C), proving
Theorem \ref{main}.
\endpf

\begin{prop}\label{remark}
Let $(S1)$ --- $(S5)$ hold. If $A(x,u)\equiv A(x)$, then
\refeq{Cesari}
   is equivalent to
 \begin{equation}\label{E322}
  E(x,y)=\bigcap_{\d>0}\coh E(x,B_\d(y)),
\end{equation}
where
$$
  E(x,y)=\set{(\z,\z^0)\in \IR\times \IR \Big| \z=f(x,y,u),\, \z^0\geq f^0(x,y,u), u\in U
  }.
  $$
\end{prop}
\textbf{Proof.}
Denote
  $$
  \tiE(x,y)=\set{(\z,\z^0)\in \IR\times \IR \Big| \z=\int_Zf(x,y,u(z))\,dz,\, \z^0\geq \int_Z f^0(x,y,u(z))\,dz, u(\cd)\in \cU_Z }.
  $$
When $A(x,v)\equiv A(x)$, \refeq{Cesari} is equivalent to
\begin{equation}
  E(x,y)=\bigcap_{\d>0}\ol{\tiE(x,B_\d(y))}.
\end{equation}
To prove \refeq{E322}, we need  only to show that
\begin{equation}
  \coh E(x,y)=\overline{\tiE(x,y)}.
\end{equation}
\textbf{I.} We first prove $\coh E(x,y)\subseteq
\overline{\tiE(x,y)}$.

For any $(\z,\z^0)\in \co E(x,y)$, there exist $\a_i$, $i=1,2,\cds,
m$ such that $\ds{\sum^m_{i=1}\a_i=1}$ and
$$
\z=\sum^m_{i=1}\a_i f(x,y,u_i),\qq\z^0\geq\sum^m_{i=1}\a_i
f^0(x,y,u_i).
$$
Define  $u_0(\cd)\in \cU_z$ by
\begin{equation}
  u_0(z)=\left\{\begin{array}{ll}
  \ds u_1, & z_1\in [0,\a_1],\\
  \ds u_k,&  z_1\in \Big(\sum^{k-1}_{i=1}\a_i,\sum^k_{i=1}\a_i\Big].
  \end{array}
  \right.
\end{equation}
Thus
$$
\int_Z f(x,y,u_0(z))\,dz=\sum^m_{i=1}\a_i f(x,y,u_i)=\z
$$
and
$$
\int _Zf^0(x,y,u_0(z))\,dz=\sum^m_{i=1}\a_i f^0(x,y,u_i)\leq \z^0.
$$
This means $(\z,\z^0)\in \tiE(x,y)$. Thus  $\co E(x,y)\subseteq
\tiE(x,y)$,  and then  $\coh E(x,y)\subseteq \overline{\tiE(x,y)}$.

\noindent\textbf{II.}  Now, we turn to prove
$\overline{\tiE(x,y)}\subseteq \coh E(x,y)$.

Let $\set{U_i^k}_{1\leq i \leq k}$ be a family of measurable
decompositions of $U$, such that

\thb{a} if $i\ne j$,  then $ \ds U_i^k\bigcap U_j^k=\emptyset $;

\thb{b} for any $k$, $\ds\bigcup^k_{j=1} U_j^k=U$;

\thb{c} $\ds \lim_{k\to +\infty}\max_{1\leq j\leq k}\diam (U_j^k)=
0$.

Moreover, let $u_i^k\in U_i^k$ and
$$
Z_i^k=\set{z\in Z \Big| u(z)\in U_i^j}.
$$
 Then by the continuity of $f$ and the lower semi-continuity of $f^0$, for a.e. $(x,y)\in \O \times \IR$,
$$
\lim_{k\to
\infty}\sum^k_{i=1}f(x,y,u_i^k)\chi_{Z_i^k}(z)=f(x,y,u(z)),\q \eqae
z\in Z
$$
and
$$
\lim_{k\to \infty}\sum^k_{i=1}f^0(x,y,u_i^k)\chi_{Z_i^k}(z)\geq
f^0(x,y,u(z)),\q \eqae z\in Z,
$$
which means for a.e. $(x,y)\in \O \times \IR$,
$$
\lim_{k\to \infty}\sum^k_{i=1}f(x,y,u_i^k)m(Z_i^k)=\int_Z
f(x,y,u(z))\,dz,
$$
and
$$
\lim_{k\to \infty}\sum^k_{i=1}f^0(x,y,u_i^k)m(Z_i^k)\geq\int_Z
f^0(x,y,u(z))\,dz.
$$
Noting that $m(E_i^k)\geq 0$ and $\ds{\sum^k_{i=1}m(Z_i^k)=1}$ for
$k=1,2,\cds$, we deduce $\tiE(x,y)\subseteq \coh E(x,y)$.
Consequently, $\overline{\tiE(x,y)}\subseteq \coh E(x,y)$. This end
the proof.
\endpf

Similar to the classical cases (see for example, Chapter 3,
Proposition 4.3 in \cite{Li}), we have the following proposition
\begin{prop}
Assume that

$(S6)$ For almost all $x\in \O$,  $f(x,\cd, v)$ is continuous
uniformly in $v\in U$ and $f^0(x,\cd, v)$ is lower semi-continuous
uniformly in  $v\in U$, i.e. for any $y\in\IR$ and $\ve>0$, there
exists a $\t=\t(x,y)>0$, such that for any $\tiy\in B_\t(y)$,
\begin{equation}\label{s6}
  \left\{
  \begin{array}{l}
    \left|f(x,\tiy,v)-f(x,y,v)\right|<\ve,\\
    f^0(x,\tiy,v)>f^0(x,y,v)-\ve,
  \end{array}
  \qq\all v\in U. \right.
\end{equation}
Then \refeq{Cesari} is equivalent to
\begin{equation}\label{GEE}
  \cE(x,y)=\overline{G\cE(x,y)},\qq\eqae\,(x,y)\in \O\times \IR.
\end{equation}
\end{prop}
\textbf{Proof.} We will prove that
\begin{equation}
  \overline{G\cE(x,y)}=\bigcap_{\d>0}\overline{G\cE(x,B_\d(y))}.
\end{equation}
In fact, we need only to show
\begin{equation*}
 \bigcap_{\d>0}\overline{G\cE(x,B_\d(y))}\subseteq \overline{G\cE(x,y)}
\end{equation*}
since
\begin{equation*}
  \overline{G\cE(x,y)}\subseteq\bigcap_{\d>0}\overline{G\cE(x,B_\d(y))}
\end{equation*}
holds obviously.

By (S6), for any $y\in \IR^n$ and $\ve>0$, there exists a
$\t=\t(x,y)>0$, such that for any $\tiy\in B_\t(y)$, \refeq{s6}
holds. For any $\d\in (0,\t)$. Let  $(P_\d,\z_\d, \z_\d^0)\in
G\cE(x,B_\d(y))$. That is
\begin{equation}
  \left\{
  \begin{array}{ll}
    \ds (P_\d)_{ij} &=\int_Z A(x,u^\d(z))(e_i+\na w^i_\d (z;x))\cd e_j\,dz,\\
  \ds  \z_\d &=\int_Z f(x,y^\d,u^\d(z))\,dz,\\
  \ds \zeta^0_\d &\geq \int_Z f^0(x,y^\d,u^\d(z))\,dz
  \end{array}
  \right.
\end{equation}
for some $u^\d(\cd)\in \cU_Z$ and $y^\d\in B_\d(y)$, where
$w^i(\cd;x)\in H^1_\#(Z)$ solves
\begin{equation*}
  \na_z\cd\Big(A(x,u^\d(z))(e_i+\na w^i_\d (z;x))\Big)=0.
\end{equation*}
Thus by \refeq{s6}, we obtain
\begin{equation}
  \left\{
  \begin{array}{l}
   \ds |\z_\d-\int_Z f(x,y,u^\d(z))\,dz|\leq \int_Z |f(x,y^\d,u^\d(z))-f(x,y,u^\d(z))|\,dz<\ve,\\
   \ds \zeta^0_\d\geq \int_Z f^0(x,y^\d,u^\d(z))\,dz>\int_Z f^0(x,y,u^\d(z))\,dz-\ve.
  \end{array}
  \right.
\end{equation}
That is,  $(P_\d,\z_\d, \z_\d^0)\in B_\ve(G\cE(x,y))$. Consequently,
\begin{equation}
  G\cE(x,B_\d(y))\subseteq B_\ve(G\cE(x,y)).
\end{equation}
Therefore
\begin{equation}
  \bigcap_{\d>0}\ol{G\cE(x,B_\d(y))}\subseteq \bigcap_{\ve>0}\ol{B_\ve(G\cE(x,y))}=\ol{G\cE(x,y)},
\end{equation}
which ends the proof.
\endpf

\end{document}